\documentclass[fleqn,11pt]{article}
\usepackage{amssymb}
\usepackage{amsmath}
\usepackage{amsthm}
\usepackage{amsopn}
\usepackage{mathrsfs}
\usepackage{harvard}
\usepackage{comment}

\numberwithin{equation}{section}

\swapnumbers

\newtheorem{satz}{Satz}[section]

\newtheorem{theorem}[satz]{Theorem}
\newtheorem{proposition}[satz]{Proposition}
\newtheorem{corollary}[satz]{Corollary}
\newtheorem{lemma}[satz]{Lemma}
\newtheorem{assumption}[satz]{Assumption}

\newtheorem{definition}[satz]{Definition}

\newtheorem{remark}[satz]{Remark}

\newtheorem{example}[satz]{Example}

\DeclareMathOperator{\E}{{\mathbb E}}

\DeclareMathOperator{\R}{{\mathbb R}}

\DeclareMathOperator{\Z}{{\mathbb Z}}
\DeclareMathOperator{\N}{{\mathbb N}}

\DeclareMathOperator{\PP}{{\mathbb P}}

\DeclareMathOperator{\spann}{span}

\DeclareMathOperator{\supp}{supp} 
 
 \DeclareMathOperator{\ran}{ran}

 \DeclareMathOperator{\KL}{KL}
 
\DeclareMathOperator{\argmin}{argmin}

\DeclareMathOperator{\Var}{Var}

\renewcommand{\phi}{\varphi}
\renewcommand{\theta}{\vartheta}
\renewcommand{\subset}{\subseteq}

\renewcommand{\cdot}{{\scriptstyle \bullet} }
\providecommand{\abs}[1]{\lvert #1 \rvert}
\providecommand{\norm}[1]{\lVert #1 \rVert}

\providecommand{\scapro}[2]{\langle #1,#2 \rangle}
\providecommand{\floor}[1]{\lfloor #1 \rfloor}

\renewcommand{\le}{\leqslant}
\renewcommand{\ge}{\geqslant}

\renewcommand{\cal}{\mathscr}     % More beautiful calligraphic

\begin{document}

\title{On rate optimality for ill-posed inverse problems\\ in econometrics\footnote{We thank Joel Horowitz,
Whitney Newey, Demian Pouzo, Yixiao Sun and the participants of the
March 2007 Oberwolfach Workshop on Semiparametrics for very helpful
discussions. X. Chen acknowledges support from the NSF/USA. The
usual disclaimer applies.} }
\author{Xiaohong Chen\footnote{Department of Economics, Yale University, Box 208281, New Haven,
CT 06520, USA. Email: xiaohong.chen@yale.edu} \and Markus
Rei{\ss}\footnote{Institute of Applied Mathematics, University of
Heidelberg, Germany. Email: reiss@statlab.uni-heidelberg.de}}
\date{First version: March 15, 2007; this version: September 10, 2007.}
\maketitle

\begin{abstract}
In this paper, we clarify the relations between the existing sets of
regularity conditions for convergence rates of nonparametric
indirect regression (NPIR) and nonparametric instrumental variables
(NPIV) regression models. We establish minimax risk lower bounds in
mean integrated squared error loss for the NPIR and the NPIV models
under two basic regularity conditions that allow for both mildly
ill-posed and severely ill-posed cases. We show that both a simple
projection estimator for the NPIR model, and a sieve minimum
distance estimator for the NPIV model, can achieve the minimax risk
lower bounds, and are rate-optimal uniformly over a large class of
structure functions, allowing for mildly ill-posed and severely
ill-posed cases.
\end{abstract}

\noindent KEY WORDS: Nonparametric instrumental regression;
Nonparametric indirect regression; Statistical ill-posed inverse
problems; Minimax risk lower bound; Optimal rate.

\setcounter{page}{0} \thispagestyle{empty}\newpage

\baselineskip=18pt

\section{Introduction}

Recently there is a growing interest in estimation for nonparametric
instrumental variables (NPIV) regression models, see e.g., Newey and
Powell (2003), Darolles, Florens and Renault (2002), Hall and
Horowitz (2005), Blundell, Chen and Kristensen (2007), Gagliardini
and Scaillet (2006), to name only a few. The estimators proposed in
these papers belong to three broad classes: (1) the finite
dimensional sieve minimum distance estimator (Newey and Powell
(2003), Ai and Chen (2003) and Blundell, Chen and Kristensen
(2007)); (2) the infinite dimensional kernel based Tikhonov
regularized estimator (Darolles, Florens and Renault (2002), Hall
and Horowitz (2005), Gagliardini and Scaillet (2006)); and (3) the
finite dimensional orthogonal series Tikhonov regularized estimator
(Hall and Horowitz (2005)). Each of these papers presents different
sets of sufficient conditions for consistency and convergence rates
of its proposed estimators. In addition, for the mildly ill-posed
case (when the singular values associated with the conditional
expectation operator decay to zero at a polynomial rate), Hall and
Horowitz (2005) establish the minimax risk lower bound in mean
integrated squared error loss for the NPIV regression model under a
set of regularity conditions that are related to their estimation
procedures. They also show that their proposed estimators achieve
this lower bound; hence their rate is optimal for the class of
structure functions they consider.

To the best of our knowledge, there is no published work that
discuss the relations among the different sets of sufficient
conditions imposed in these various papers. Therefore, it is unclear
whether the minimax risk lower bound derived in Hall and Horowitz
(2005) is still the lower bound under regularity conditions stated
in the other papers. It is also unclear whether the estimators
proposed in the other papers are rate optimal in a minimax framework
corresponding to the conditions stated in these papers. Moreover,
when the NPIV problem is severely ill-posed (for instance, when the
singular value associated with the conditional expectation operator
decays to zero at an exponential rate), there are no published
results on minimax rates.

In this paper, we address these issues based on a general
formulation of the problems. In Section 2, we first present the NPIV
models. We then provide two basic regularity conditions: the
approximation and the link conditions. The approximation condition
is about the complexity of the class of the structural functions,
which is measured as the best finite dimensional linear
approximation error rate in terms of a basis expansion that may not
be the eigenfunction basis of the conditional expectation operator.
The link condition is about the relative smoothness of the
conditional expectation operator in terms of the basis used in the
first condition. We show that these two regularity conditions are
natural generalizations of, and are automatically satisfied by, the
so-called ``general source condition'', an assumption commonly
imposed in the literature on ill-posed inverse problems. Our two
basic regularity conditions are also implied by the ones assumed in
the literature on NPIV models, such as those imposed in Darolles,
Florens and Renault (2002), Hall and Horowitz (2005), and Blundell,
Chen and Kristensen (2007). In Section 3, we first show that the
NPIV model is no more informative than the reduced form
nonparametric indirect regression (NPIR) model (actually the model
assuming a known conditional expectation operator of the endogenous
regressor given the instrumental variables). Under the two basic
regularity conditions stated in Section 2, we derive the minimax
risk lower bound in mean integrated squared error loss for the NPIR
and the NPIV models, allowing for both the mildly ill-posed case and
the severely ill-posed case. In Section 4, we present a simple
projection estimator for the NPIR models, and establish that it
achieves the lower bounds and hence is rate-optimal in the minimax
sense. When restricting our conditions to various special cases,
including the nonparametric mean regression models and the NPIR
models under general source conditions, our results reproduce the
existing known minimax optimal rates for these special cases. But
more importantly, our minimax optimal rate results cover many new
cases as long as their model specifications satisfy the
approximation and the link conditions. We also discuss what could
happen if the link condition on the relative smoothness of the
conditional expectation operator is not satisfied. In Section 5, we
show that the sieve minimum distance (SMD) estimator for the NPIV
models is rate-optimal in the minimax sense. In fact, we show that
both the projection estimator for the NPIR models and the SMD
estimator for the NPIV models are rate-optimal uniformly over a
large class of structure functions, allowing for arbitrarily
decaying speed of the singular values of the conditional expectation
operator. Section 6 provides some further discussions on the
regularity conditions. Section 7 briefly concludes, and all the
proofs are gathered in the Appendix.

Before we conclude this introduction, we mention closely related
work in more abstract settings of linear ill-posed inverse problems.
First, there exist many papers and some monographs devoted to
constructing estimators and deriving optimal convergence rates in
the deterministic noise framework with a known operator (or a known
operator up to a deterministically perturbed error with a specified
error rate). See, e.g., Engl, Hanke and Neubauer (1996), Nair,
Pereverzev and Tautenhahn (2005) and the references therein. Second,
there are also many results on minimax optimal rates in mean
integrated squared error loss in the random white noise framework
with a known operator; see, e.g., Cohen, Hoffmann and Rei\ss\
(2004), Bissantz, Hohage, Munk and Ruymgaart (2007) and the
references therein. Third, there are a few recent papers on
constructing estimators that achieve optimal convergence rates in
the presence of a white noise and an unknown operator, but assuming
the existence of an estimator of the operator with a rate. See,
e.g., Efromovich and Koltchinskii (2001) and Hoffmann and Rei\ss\
(2007).

\section{NPIV models and basic regularity conditions}

We first specify the NPIV regression model as
\begin{equation}
\label{npiv}
Y_i=h_0(X_i)+U_i,\quad \E[U_i\,|\,W_i]=0, \quad i=1,\ldots,n,
\end{equation}%
with observations $\{(X_i,Y_i,W_i)\}_{i=1}^{n}$, a random sample
from the unknown joint distribution of $(X,Y,W)$. Here $Y$ is a
scalar dependent variable, $X$ is a vector of endogenous regressors
in ${\R}^d$ and $W$ is a vector of instrumental variables in
${\R}^d$ that satisfy the property $\E[U\,|\,W ]=0$. (For the ease
of presentation we assume that $X$ and $W$ do not contain any common
variables, and the conditional density of $X$ given $W$ is
well-defined). The parameter of interest is the unknown structure
function $h_0(\cdot)$, while the joint law ${\cal L}_{UWX}$ of
$(U,W,X)$ is an unknown nuisance function.

Let us introduce the Hilbert spaces
\begin{align*}\textstyle
L^2_X&=\{h:{\R}^d \to\R\,|\,\norm{h}_X^2:=\E[h^2(X)]<\infty\},\\
L^2_W&=\{g:{\R}^d \to\R\,|\,\norm{g}_W^2:=\E[g^2(W)]<\infty\}.
\end{align*}
Since the conditional distribution of $X$ given $W$ is unspecified,
the conditional expectation operator
\[(Kh)(w):=\E[h(X)\,|\,W=w]\]
is unknown, except that it is an integral operator mapping from
$L^2_X$ to $L^2_W$. This operator is the key in the construction of
estimators of $h_0$ because by conditioning on $W$ in \eqref{npiv}
and using $\E[U\,|\,W]=0$ we obtain
\[ \E[Y\,|\,W]=\E[h_0(X)\,|\,W]+\E[U\,|\,W]=Kh_0(W).\]
Consequently, by regressing $Y$ on $W$, estimating $K$ and using
this relationship we can hope to retrieve an estimator of $h_0$.

Let ${\cal H}$ denote a subset of $L^2_X$ and assume $h_0 \in {\cal
H}$. Here ${\cal H}$ captures all the prior information (such as the
smoothness and/or shape properties) about the unknown structure
function $h_0$. To ensure that there is a unique solution $h_0 \in
{\cal H}$ for the NPIV model (\ref{npiv}), in this paper we assume
that the operator $K$ satisfies the following restriction:
\begin{equation} \label{identification}
 \{h\in {\cal H} : Kh=0\}=\{0\}.
\end{equation}%
Depending on the choice of the function class $\cal H$, the
identification condition (\ref{identification}) imposes different
restrictions on the operator $K$ (or equivalently, on the
conditional density of $X$ given $W$). For example, if $\cal H =
L^2_X$, then condition (\ref{identification}) becomes the standard
identification condition that $K$ is \emph{injective}, i.e., ${\cal
N}(K):=\{h\in L^2_X : Kh=0\}=\{0\}$, (or equivalently, the
conditional density of $X$ given $W$ is complete); see, e.g., Newey
and Powell (2003), Darolles, Florens and Renault (2002), Carrasco,
Florens and Renault (2007). If $\cal H =\{h\in L^2_X: \sup_x
|h(x)|\le 1\}$, then condition (\ref{identification}) corresponds to
assume that the conditional density of $X$ given $W$ is bounded
complete; see, e.g., Chernozhukov and Hansen (2005), Chernozhukov,
Imbens and Newey (2007), Blundell, Chen and Kristensen (2007). For
additional results on identification in semi/nonparametric models
with endogeneity, see, e.g., Blundell and Powell (2003), Florens
(2003), Florens, Johannes and Van Bellegem (2007) and the references
therein.

\subsection{Basic regularity conditions}

In this paper we would like to establish a minimax risk lower bound
for the NPIV model, that is, we would like to derive a result of the
form: there are a finite constant $c>0$ and a rate function
$\delta_n \downarrow 0$ as $n\uparrow\infty$ such that
\[\lim_{n\to\infty} \Big(\delta_n^{-1}\inf_{\hat h_n} \sup_{h\in{\cal H}} \E_{({\cal
L}_{UWX},h)} [\norm{\hat h_n-h}_X^2]\Big) \ge c
\]
where the infimum is over all possible estimators $\hat h_n$ for $h
\in {\cal H}$. Note that a NPIV model (\ref{npiv}) is completely
specified by prescribing the joint law ${\cal L}_{UWX}$ of $(U,W,X)$
and the structure function $h$. This lower bound $\delta_n$ will be
valid for quite general forms of ${\cal L}_{UWX}$, independently of
knowing or not knowing it. In particular, although the mean squared
error loss and the class of structure functions ${\cal H}$ will be
defined in terms of the distribution of $X$, there is no need to
assume any explicit properties of this distribution to derive a
minimax lower bound.

We would also like to present some particular estimators that attain
the lower bound rate $\delta_n$. However, before we could establish
any minimax lower and upper bounds, it is clear that we have to
impose some conditions on the class of structure functions $\cal H$
and on the conditional expectation operator $K$. In this paper, we
implicitly assume that the prior information about $\cal H$ already
includes some regularity properties that could be described in terms
of a Hilbert scale generated by a conveniently chosen (by the
researcher) operator $B$. The regularizing action of the conditional
expectation operator $K$ would also be described as some smoothness
relative to the known operator $B$. Formally, let $B: Dom(B)\subset
L_X^2\to L_X^2$ be a densely defined self-adjoint, strictly positive
definite, and \emph{unbounded} operator (such as differential
operators with boundary constraints). For the ease of presentation
we assume that $B$ has eigenvalues $\nu_k\uparrow\infty$ with
corresponding $L_X^2$-normalized eigenfunctions $\{u_k\}$ which then
form an orthonormal basis of $L_X^2$. For non-discrete spectrum our
results will still hold, but the presentation would become more
technical, using spectral measures and abstract functional calculus.

Throughout this paper we denote by ${\cal H}(r,R)$ a subset of $\cal
H \subset L_X^2$, and assume the following:
\begin{assumption}[approximation condition]
\label{h-smooth} There are finite constants $r,R>0$ such that ${\cal
H}(r,R)$ consists of functions $h$ satisfying
\begin{equation} \label{h-basis}
\inf_{\{a_k\}:\sum_{k} a_k^2 <\infty} \norm{h-\sum_{k=1}^m a_k
u_k}_X^2 \le R^2 \nu_{m+1}^{-2r} \quad \text{for all }~m\in\N.
\end{equation}%
\end{assumption}

Note that the left-hand side of \eqref{h-basis} gives the error in
approximating $h$ optimally by an element of the $m$-dimensional
space spanned by the basis functions $\{u_1,\ldots,u_m\}$. So,
Assumption \ref{h-smooth} characterizes the regularity (or
smoothness) of the structure functions in ${\cal H}(r,R)$ by the
$L^2_X$-approximation error rates when they are approximated by the
basis $\{u_k\}$ associated with $B$. Clearly, Assumption
\ref{h-smooth} will give a bound on the bias and implies that ${\cal
H}(r,R)$ is a compact set in $L^2_X$. For many typical smooth
function classes and basis functions like the Fourier basis,
wavelets or splines the approximation error rates are well known.
%In particular, this assumption is implied by ${\cal H}(r,R)$ being an
%$L^2$-Sobolov ball or a H\"older ball of regularity $r$.

For any $s>0$ and $h\in Dom(B^s)\subset L_X^2$ we write
$\norm{h}_s:=\norm{B^s h}_X$. Let $H^s$ denote the completion of
$Dom(B^s)$ under the norm $\norm{\cdot}_s$. $\{H^s\}_{s>0}$ is
called a {\em Hilbert scale} generated by $B$ (see, e.g., Engl,
Hanke and Neubauer (1996) for its detailed properties). For any
finite constants $r, R>0$, we define a Sobolev-type ellipsoid as
$H_R^r :=\{h\in H^r ,~ \norm{h}_r\le R\}$. Since
\[ H_R^r
=\Big\{h=\sum_{k=1}^{\infty} \scapro{h}{u_k}_X u_k ,~ \norm{h}_r^2 =
\sum_{k=1}^{\infty}\nu_k^{2r} \scapro{h}{u_k}_X^2 \le R^2 \Big\},
\]
it is clear that $H_R^r$ is a subset of ${\cal H}(r,R)$. It is also
easy to see that the following hyperrectangle ${\Theta}^r_{R'}$ in
$L^2_X$ is a subset of ${\cal H}(r,R)$ for $R'>0$ sufficiently
small:
\[
{\Theta}^r_{R'} :=\Big\{h=\sum_{k=1}^{\infty} \scapro{h}{u_k}_X u_k
,~ |\scapro{h}{u_k}_X|\le R' \nu_k^{-\beta}\Big\},\quad \beta
=r+\tfrac{1}{2} > \tfrac{1}{2}.
\]
Let us now formulate the mapping properties of the conditional
expectation operator $K$ in terms of the (generalized) Hilbert scale
generated by $B$.

\begin{assumption}[link condition]
\label{K-smooth} There are a continuous increasing function
$\phi:\R^+\to\R^+$ and a constant $M>0$ such that $\norm{Kh}_W\le
M\norm{[\phi(B^{-2})]^{1/2}h}_X$ for all $h\in L_X^2$.
\end{assumption}
Assumption \ref{K-smooth} is in fact equivalent to the range
inclusion condition:
\[\ran(\abs{K})\subset
\ran([\phi(B^{-2})]^{1/2})\text{ with } \abs{K} : =(K^\ast K)^{1/2},
\]
where $K^\ast$ denotes the Hilbert space ($L^2_X$) adjoint of $K$.
For the NPIV models, under mild conditions, the conditional
expectation operator $K$ is a compact operator. Thus the
self-adjoint compact operator $K^\ast K$ has the
eigenvalue-eigenfunction decomposition $\{\lambda_k, e_k\}$, where
the eigenvalues are arranged in non-increasing order: $\lambda_k \ge
\lambda_{k+1} \ge ...
>0$ and $\lambda_k $ tends to zero as $k\uparrow\infty$.
Then Assumption \ref{K-smooth} can be equivalently restated in terms
of two possibly different orthonormal bases $\{e_k\}$ and $\{u_k\}$
of $L^2_X$:
\begin{equation} \label{K-basis}
\sum_{k=1}^{\infty} \lambda_k \scapro{h}{e_k}_X^2 \le
M^2\sum_{k=1}^{\infty} \phi(\nu^{-2}_k ) \scapro{h}{u_k}_X^2 \quad
\text{for all }~ h\in L^2_X.
\end{equation}%

\begin{remark}
We can rewrite Assumptions \ref{h-smooth} and \ref{K-smooth} without
specifying the operator $B$ explicitly. All we require are the
existence of an orthonormal basis $\{u_k\}$ in $L^2_X$ and a
sequence of increasing positive real numbers $\{\nu_k\}$ such that
equations (\ref{h-basis}) and (\ref{K-basis}) hold. In fact, we can
then construct the self-adjoint unbounded operator $B$ according to
\[ Bh=\sum_{k=1}^\infty \nu_k \scapro{h}{u_k}_X u_k ,\]
with $Dom(B)=\{h\in L^2_X : \sum_{k=1}^\infty \nu_k^2
\scapro{h}{u_k}_X^2 < \infty\}$.
\end{remark}

\begin{example}
Suppose that $X$ is uniformly distributed on the interval $[0,1]$
and let $Bf(x):=-f''(x)$ for all $f\in L^2([0,1])$ with $f''\in
L^2([0,1])$ and with periodic boundary conditions. Then $B$ has
(complex-valued) eigenfunctions $u_k(x)=\exp(2\pi kix)$ with
eigenvalues $\nu_k=(2\pi k)^2$ such that
\[ H^r=\{f\in L^2([0,1])\,:\,
\sum_{k\in\Z}\nu_k^{2r}\abs{\scapro{f}{u_k}}^2<\infty\}
\]
is the classical $L^2$-Sobolev space $H^{2r}_{per}$ of regularity
(smoothness) $2r$ with periodic boundary conditions. See, e.g.,
Edmunds and Evans (1987) for many examples of generating smooth
function spaces from differential operators.

For the typical choice $\phi(t)=t^a$ for some $a>0$, Assumption
\ref{K-smooth} translates to $\norm{Kh}_W\le M\norm{B^{-a}h}_X$,
which means intuitively that the operator $K$ regularizes at least
as much as $B^{-a}$. In the case $Bf(x):=-f''(x)$ the operator $K$
acts like integrating at least $(2a)$-times, i.e. maps $L^2$ to the
$L^2$-Sobolev space of regularity $2a$.
\end{example}

In the statistics literature, for the standard nonparametric mean
regression model (i.e., the model in which $K$ is the identity
operator), the minimax risk lower and upper bounds have been
established in mean integrated squared error loss for various
classes of functions $\cal H$ such as a Sobolev ball (ellipsoid), a
H\"older ball (hyperrectangle) or a Besov ball (ellipsoid or
hyperrectangle or Besov body); see, e.g., Donoho, Liu and MacGibbon
(1990), Yang and Barron (1999) and the references therein. As shown
in these papers, what matters for minimax risk lower and upper
bounds for nonparametric mean regression estimation is the
complexity of the class of functions $\cal H$ that can be measured
in terms of best finite dimensional approximation numbers. This
motivates us to impose Assumption \ref{h-smooth}. However, since the
basis $\{u_k\}$ (of the operator $B$) used to construct the best
finite dimensional approximations for the class of functions $\cal
H$ may differ from the eigenfunction basis $\{e_k\}$ (of the
operator $K^\ast K$), we have to impose Assumption \ref{K-smooth} to
link these two.

We shall refer to Assumptions \ref{h-smooth} and \ref{K-smooth} as
the two basis regularity conditions; and sometimes call Assumption
\ref{h-smooth} the approximation condition and Assumption
\ref{K-smooth} the link condition. Both assumptions are satisfied by
the ones imposed in the literature, such as those in Cohen, Hoffmann
and Rei\ss\ (2004), Efromovich and Koltchinskii (2001), Hoffmann and
Rei\ss\ (2007), Blundell, Chen and Kristensen (2007), Chen and Pouzo
(2007) and others. In the next subsection we show that these two
basic regularity conditions are automatically satisfied by the
so-called ``general source condition'', which in turn are satisfied
by conditions imposed in Hall and Horowitz (2005) and all the other
papers using the general source condition.

\subsection{Relation to source conditions}

In the numerical analysis literature on ill-posed inverse problems
it is common to measure the smoothness (regularity) of the function
class $\cal H$ according to the spectral representation of the
operator $K^{\ast}K$. Denote by $\norm{K}:=\sup_{h: \norm{h}_X \le
1}{\norm{Kh}_W}$ the operator norm. The so-called ``\emph{general
source condition}'' assumes that there is a continuous function
$\psi$ defined on $[0,\norm{K}^2]$ with $\psi (0)=0$ and
${\lambda}^{-1/2}\psi(\lambda)$ non-decreasing such that
\begin{equation} \label{source}
{\cal H}_{source} :=\Big\{h=\psi (K^\ast K) g, ~g\in L^2_X,
~\norm{g}_X^2 \le R\Big\},~\text{for a finite constant} ~R,
\end{equation}%
and the original ``source condition'' corresponds to the choice
$\psi (\lambda)= {\lambda}^{1/2}$ (see Engl, Hanke and Neubauer
(1996)). If $K^\ast K$ is compact with eigenvalue-eigenfunction
system $\{\lambda_k, e_k\}$, then (\ref{source}) is equivalent to
\[
{\cal H}_{source} =\Big\{h=\sum_{k=1}^{\infty} \scapro{h}{e_k}_X e_k
,~ \sum_{k=1}^{\infty}\frac{\scapro{h}{e_k}_X^2}{\psi^2 (\lambda_k)}
\le R^2\Big\}.
\]
Therefore the general source condition implies our Assumptions
\ref{h-smooth} and \ref{K-smooth} by setting $u_k = e_k$ and
$\nu_k^{-r} = \psi (\lambda_k)$ for all $k\ge 1$, and
$\phi(B^{-2})=K^\ast K$.

In the econometrics literature on NPIV estimation, Darolles, Florens
and Renault (2002) impose a smoothness condition on the true
structure function $h_0$ that is closely related to the source
condition. Precisely, they assume $h_0 \in {\cal H}_{DFR}$, where
\begin{equation} \label{DFR}
{\cal H}_{DFR} =\Big\{h \in L^2_X,~
\sum_{k=1}^{\infty}\frac{\scapro{h}{e_k}_X^2}{(\lambda_k)^{a}}
<\infty \Big\},\quad \text{for some} ~a \ge 1.
\end{equation}%
Darolles, Florens and Renault (2002) use this assumption $h_0 \in
{\cal H}_{DFR}$ to establish the convergence rate of their
kernel-based Tikhonov regularized estimator in mean squared error
metric $\E_{h_0}[\norm{\hat{h}-h_0}_X^2]$. This rate, however, will
not hold uniformly over $h_0 \in {\cal H}_{DFR}$, since the series
in \eqref{DFR} is not uniformly bounded away from infinity, which is
the role of $R \in (0,\infty )$ in the definition of ${\cal
H}_{source}$.

Hall and Horowitz (2005) assume that $h_0$ belongs to a
hyperrectangle in $L^2_X$, using the eigenfunctions $\{e_k\}$ of the
operator $K^\ast K$ as a basis:
\begin{equation} \label{HH}
{\cal H}_{HH} =\Big\{h=\sum_{k=1}^{\infty} \scapro{h}{e_k}_X e_k ,~
|\scapro{h}{e_k}_X|\le R' k^{-\beta}\Big\},
\end{equation}%
which, when $\beta>1/2$ plays the role of $r+1/2$, implies our
Assumptions \ref{h-smooth} and \ref{K-smooth} by setting $u_k =
e_k$, $\nu_k = k$ for all $k\ge 1$, and $\phi(B^{-2})=K^\ast K$. In
addition, Hall and Horowitz (2005) also assume that the eigenvalues
$\{\lambda_k\}$ of the operator $K^\ast K$ are such that $\lambda_k
\ge const. k^{-\alpha}$ for some $\alpha
> 1$ and $2\beta > \alpha \ge \beta -\tfrac{1}{2}$, which suggests that
we could set $\phi (t)=t^{\alpha /2}$.

%For the linear ill-posed
%inverse problem with a known operator $K$ and when the problem is
%mildly ill-posed, K (1993) derives minimax rate assuming $\cal H$ is
%a Sobolev ellipsoid.

\section{The lower bound}

In this section we shall establish a minimax risk lower bound for
the NPIV model under the two basic regularity conditions stated in
Section 2. We derive this result by first establishing that the NPIV
model is no more informative than the reduced form nonparametric
indirect regression (NPIR) model. First, the following abstract
assumption ensures a certain complexity of the statistical NPIV
model and permits the residuals of $Y$ given $W$ to be Gaussian.
Recall that $\cal L_Z$ denotes the law of the random vector $Z$.

\begin{assumption}\label{ass-complexity}
Let $\sigma_0>0$ be a finite constant. Let $\cal C$ be a (possibly
very large) set of elements $({\cal L}_{UWX},h)$ such that the
following property holds:
\begin{itemize}
\item For all $h\in {\cal H}$, there is a law ${\cal L}_{UWX}$ with $({\cal L}_{UWX},h)\in{\cal C}$
such that ${\cal L}_{WY}$ is determined by ${\cal L}_{UWX}$ and $h$,
and that
\[V_i:=Y_i-\E[Y_i\,|\,W_i]=h(X_i)-(Kh)(W_i)+U_i
\]
given $W_i$ is $N(0,\sigma^2(W_i))$-distributed with
$\sigma^2(W_i)\ge \sigma_0^2$.
\end{itemize}
\end{assumption}

\begin{example}
A typical NPIV model (\ref{npiv}) satisfying Assumption
\ref{ass-complexity} is generated by taking $W_i$ from an arbitrary
probability distribution ${\cal L}_{W}$, then generating $X_i$
according to a conditional density of $X$ given $W$, generating
$V_i$ according to $N(0,\sigma^2(W_i))$, and defining
\[ U_i:=(Kh)(W_i)-h(X_i)+V_i.\]
\end{example}

\subsection{Reduction from NPIV model to NPIR model}

For each NPIV model, we specify the reduced form NPIR model as
\[ Y_i=(Kh)(W_i)+V_i,\quad i=1,\ldots,n,\]
with $(W_i,V_i)$ i.i.d., ${\cal L}_{V|W=w}=N(0,\sigma^2(w))$, $h \in
\cal H$ the unknown structure function, and $K$ a {\em known}
injective operator from $L^2_X$ to $L^2_W$. The observations
corresponding to the NPIR are $\{(Y_i,W_i)\}_{i=1}^{n}$. We shall
now formally prove,
 that the NPIV model is statistically more demanding than an
indirect regression model with known operator $K$. We compare
statistical experiments in a decision-theoretic sense (see Le Cam
and Yang (2000)), and therefore, have to ensure first that the
classes of parameters are compatible.

\begin{definition}
Let Assumption \ref{ass-complexity} hold. The NPIR model class
${\cal C}_0$ consists of all model parameters $({\cal
L}_{W'},\sigma(\cdot),h)$ such that there is $({\cal L}_{UWX},h)\in
{\cal C}$ with the following properties:  ${\cal L}_{W}={\cal
L}_{W'}$, $\sigma^2(w)\ge \sigma_0^2>0$, the conditional law ${\cal
L}_{X|W}$ is prescribed according to $K$, and ${\cal L}_{U|WX}$ is
arbitrary among the conditions imposed in $\cal C$.
\end{definition}

\begin{lemma}\label{LemReduction}
The NPIR model is more informative than the NPIV model in the sense
that for each estimator $\hat h_n$ for the NPIV model there is an
estimator $\tilde h_n$ for the NPIR model with
\[ \sup_{({\cal L}_W,\sigma(\cdot),h)\in{\cal C_0}} \E_{({\cal L}_W,\sigma(\cdot),h)}[\norm{\tilde h_n-h}_X^2]
\le \sup_{({\cal L}_{UWX},h)\in{\cal C}} \E_{({\cal L}_{UWX},h)}
[\norm{\hat h_n-h}_X^2].
\]
\end{lemma}

\subsection{The lower bound}

We now formally present the minimax risk lower bound for the NPIR
and the NPIV models in mean squared error loss. We establish the
lower bound by considering asymptotically least favorable Bayes
priors, more specifically, by applying Assouad's cube technique; see
e.g. Korostelev and Tsybakov (1993) or Yang and Barron (1999). In
this paper we use the notation $a_n \asymp b_n$ to mean that there
is a finite positive constant $c$ such that $c a_n \le b_n \le
c^{-1} a_n$.

Since $\sup_{h\in{\cal H}(r,R)}\E_h[\norm{\hat h_n-h}_X^2]\ge
\sup_{h\in H^r_R}\E_h[\norm{\hat h_n-h}_X^2]$, it suffices to
establish the lower bound for functions in $H_R^r$, a subset of
${\cal H}(r,R)$.

\begin{theorem}\label{thm-lb}
Let Assumptions \ref{h-smooth} and \ref{K-smooth} hold. For the NPIR
model we have the following minimax risk lower bound:
\[
\inf_{\hat h_n}\sup_{h\in H^r_R}\E_h[\norm{\hat h_n-h}_X^2]\ge
\tfrac{\sigma_0^2}{4\exp(4M)}\delta_n, \quad \delta_n :=
n^{-1}\sum_{k=1}^m [\phi(\nu_k^{-2})]^{-1},
\]
where the infimum runs over all possible estimators $\hat{h}_n$
based on $n$ observations, and $m$ is the largest possible integer
satisfying
\[\sigma_0^2n^{-1}\sum_{k=1}^m\nu_k^{2r}[\phi(\nu_k^{-2})]^{-1}\le R^2.
\]
(1) Mildly ill-posed case: Let $\phi(t)=t^a$ and $\nu_k\asymp
k^\epsilon$ for some $a, \epsilon >0$. If $m\asymp n^{1/(2r\epsilon
+2a\epsilon+1)}$, then $\delta_n \asymp
n^{-2r/(2r+2a+\epsilon^{-1})}$.

(2) Severely ill-posed case: Let $\phi(t)=\exp(-t^{-a/2})$,
$\nu_k\asymp k^\epsilon$ for some $a, \epsilon>0$. If
$m=c\log(n)^{1/a\epsilon}$ with a sufficiently small $c>0$, then
$\delta_n \asymp (\log n)^{-2r/a}$.
\end{theorem}

The next corollary follows directly from Lemma \ref{LemReduction}
and Theorem \ref{thm-lb}; hence we omit its proof.

\begin{corollary}
Let Assumptions \ref{h-smooth}, \ref{K-smooth} and
\ref{ass-complexity} hold. For the NPIV model we have the same
minimax risk lower bound:
\[
\inf_{\hat h_n}\sup_{h\in H^r_R}\E_h[\norm{\hat h_n-h}_X^2]\ge
\tfrac{\sigma_0^2}{4\exp(4M)} \delta_n,\quad \text{with}~\delta_n
~\text{given in Theorem \ref{thm-lb}},
\]
where the infimum runs over all possible estimators $\hat{h}_n$
based on $n$ observations.
\end{corollary}

\begin{remark}
For the proof of the lower bound we have to consider the likelihood
between the observations. This is why we require Gaussianity.
Nevertheless, the proof works the same for other error densities,
but bounding the Kullback-Leibler or Hellinger distance between
alternatives might be more cumbersome.

Let us also mention that the proof strategy can also yield a lower
bound for convergence in probability:
\[
\inf_{\hat h_n}\sup_{h\in H^r_R}P_h\Big(\delta_n^{-1}\norm{\hat
h_n-h}_X^2\ge \tfrac{\sigma_0^2}{4\exp(4M)}\Big)\ge c>0 ,\quad
\text{with}~\delta_n ~\text{given in Theorem \ref{thm-lb}},
\]
cf. Korostelev and Tsybakov (1993).
\end{remark}

Note that Assumption 2.2 is automatically satisfied under the
general source condition with $K^\ast K = \phi(B^{-2})$. Following
the proof of Theorem \ref{thm-lb}, we immediately obtain:
\begin{remark}\label{source-lb}
Suppose that Assumption \ref{h-smooth} is satisfied with $h\in {\cal
H}_{source}$ and $u_k = e_k$, $\nu_k^{-r} = \psi (\lambda_k)$ for
all $k\ge 1$. Let $\phi(B^{-2})=K^\ast K$. Then, for NPIR model and
for NPIV model (under Assumption \ref{ass-complexity}), we have the
same minimax risk lower bound:
\[
\inf_{\hat h_n}\sup_{h\in {\cal H}_{source}}\E_h[\norm{\hat
h_n-h}_X^2]\ge \tfrac{\sigma_0^2}{4\exp(4M)} \delta_n,\quad
\text{with}~\delta_n ~\text{given in Theorem \ref{thm-lb}},
\]
where the infimum runs over all possible estimators $\hat{h}_n$
based on $n$ observations. An equivalent way to determine the lower
bound $\delta_n$ is to choose the largest possible integer $m$ such
that
\[
\delta_n =n^{-1}\sum_{k=1}^m \lambda^{-1}_k,\quad
\sigma_0^2n^{-1}\sum_{k=1}^m [\psi (\lambda_k)]^{-2}\lambda^{-1}_k
\le R^2.
\]
\end{remark}

\section{An upper bound for the NPIR model}

We prove an upper bound for the NPIR model. The aim of this section
is to convince the reader that the lower bounds given in Section 3
are rate-optimal, and to provide an easy method to attain these
rates. Again we assume that $B$ has eigenvalues
$\nu_k\uparrow\infty$ with corresponding $L_X^2$-normalized
eigenfunctions $(u_k)$ which then form an orthonormal basis of
$L_X^2$. For $m\ge 1$ we define our estimator as
\begin{equation}
\label{npir-estimate} \hat h_n:=\sum_{k=1}^m\hat\eta_ku_k, \quad
\hat\eta_k:=\frac1n\sum_{i=1}^nY_i((K^\ast)^{-1}u_k)(W_i).
\end{equation}
This simple projection procedure using the basis $\{u_k\}$ (of $B$)
does not seem to have been studied before. It is a natural
generalization of the well-known spectral cut-off method using the
eigenfunction basis $\{e_k\}$ of $K^\ast K$. Given the prior
information about $\cal H (r,R)$, this is a mathematically
satisfactory construction.

For the upper bound we impose the following assumptions on the NPIR
model.
\begin{assumption}\label{npir-a1}\mbox{}\\
(1) There is a finite $\sigma_1>0$ such that
$\sigma(w)\le\sigma_1$ for all $w\in\supp({\cal L}_W)$;\\
(2) There is a finite $S>0$ such that
$\norm{Kh}_\infty=\sup_{w\in\supp({\cal L}_W)}\abs{(Kh)(w)}\le S$
for all $h \in \cal H (r,R)$.
\end{assumption}
Assumption \ref{npir-a1} is typically assumed in papers on
nonparametric estimation of ill-posed indirect regression; see, e.g.
Bissantz, Hohage, Munk and Ruymgaart (2007). When $K$ is the
identity operator, Assumption \ref{npir-a1}(2) becomes to require
that $\norm{h}_{\infty} \le S$ for all $h \in \cal H (r,R)$, which
is a condition imposed in Yang and Barron (1999, theorems 6 and 7)
to derive their minimax rate for a standard nonparametric regression
model.

\begin{assumption}[reverse link condition]
\label{npir-a2} There is a finite $c>0$ such that $\norm{K h}_W\ge
c\norm{[\phi(B^{-2})]^{1/2}h}_X$ for all $h\in L_X^2$.
\end{assumption}
Assumption \ref{npir-a2} is the reverse condition of Assumption
\ref{K-smooth} and is often imposed in papers on ill-posed inverse
problems. We shall sometimes call Assumptions \ref{K-smooth} and
\ref{npir-a2} together as the exact link (or exact range) condition.
See Subsection 4.2 for a relaxation of this condition.

\begin{proposition}\label{prop-npir}
For the NPIR models, suppose that Assumptions \ref{h-smooth},
\ref{npir-a1} and \ref{npir-a2} hold. Then the estimator $\hat h_n$
defined in (\ref{npir-estimate}) satisfies
\[
\sup_{h\in{\cal H}(r,R)}\E_h[\norm{\hat h_n-h}_X^2]\le
\nu_{m+1}^{-2r}R^2 + 2n^{-1}(S^2+\sigma_1^2)c^{-2}\sum_{k=1}^m
[\phi(\nu_k^{-2})]^{-1}.
\]
If $m=m(n)$ is such that
$n^{-1}\sum_{k=1}^m\nu_k^{2r}[\phi(\nu_k^{-2})]^{-1}\asymp 1$, then,
under Assumption \ref{K-smooth}, this estimator $\hat h_n$ is
rate-optimal in the minimax sense: there is a finite constant $C>0$
such that
\[
\sup_{h\in{\cal H}(r,R)}\E_h[\norm{\hat h_n-h}_X^2]\le
C\nu_{m+1}^{-2r} \asymp n^{-1}\sum_{k=1}^m [\phi(\nu_k^{-2})]^{-1}
\asymp \delta_n,\quad \text{with}~\delta_n ~\text{given in Theorem
\ref{thm-lb}}.
\]

(1) Mildly ill-posed case: Let $\phi(t)=t^a$ and $\nu_k\asymp
k^\epsilon$ for some $a, \epsilon >0$. If $m\asymp n^{1/(2r\epsilon
+2a\epsilon+1)}$, then $\delta_n \asymp
n^{-2r/(2r+2a+\epsilon^{-1})}$.

(2) Severely ill-posed case: Let $\phi(t)=\exp(-t^{-a/2})$,
$\nu_k\asymp k^\epsilon$ for some $a, \epsilon>0$. If
$m=c\log(n)^{1/a\epsilon}$ with a sufficiently small $c>0$, then
$\delta_n \asymp (\log n)^{-2r/a}$.
\end{proposition}

\begin{remark}
As the proof reveals, the upper bound does not require that the
errors are Gaussian, the existence of second moments suffices.
\end{remark}

\begin{remark}\label{identity-minimax}
When $K$ is the identity operator, the NPIR model becomes the
standard nonparametric mean regression model, and Assumption
\ref{K-smooth} is automatically satisfied with $\phi ()$ being a
constant, then Theorem \ref{thm-lb} and Proposition \ref{prop-npir}
together reproduce the well-known minimax lower and upper bounds for
the nonparametric mean regression model (see, e.g., theorem 7 of
Yang and Barron (1999)), in which $\delta_n \asymp \frac{m}{n}$, and
$m$ is the largest possible integer satisfying $\nu_{m+1}^{-2r}
\asymp \frac{m}{n}$.
\end{remark}

Comparing the minimax optimal rates in mean integrated squared error
loss for the nonparametric mean regression model and for the NPIR
model, we see the squared bias is of the same order
($\nu_{m+1}^{-2r}$), but the variance blow up from $\frac{m}{n}$ for
the nonparametric mean regression model to $n^{-1}\sum_{k=1}^m
[\phi(\nu_k^{-2})]^{-1}$ for the NPIR model.

Notice that the minimax optimal rate $\delta_n \asymp (\log
n)^{-2r/a}$ for the severely ill-posed case is independent of
$\epsilon$ (hence independent of the dimension $d$ of $X$). For the
mildly ill-posed case, when $\phi(t)=t^a$ and $\nu_k\asymp
k^\epsilon$ for some $a>0$ and $\epsilon=1/d$, Theorem \ref{thm-lb}
and Proposition \ref{prop-npir} together give the minimax optimal
rate $\delta_n=n^{-2r/(2r+2a+d)}$ for the NPIR models. This rate is
well known for the special case when $\cal H (r,R)$ is a
$d$-dimensional Sobolev ball $H^r_R$ and the operator $K$ is
elliptic with ill-posedness degree $a$ (i.e., $\norm{Kh}_W \asymp
\norm{B^{-a}h}_X$ for all $h\in L^2_X$); see, e.g., Cohen, Hoffmann
and Rei\ss\ (2004).

Note that Assumptions 2.2 and 4.2 are automatically satisfied under
the general source condition with $K^\ast K = \phi(B^{-2})$.
Applying Proposition \ref{prop-npir} and Remark \ref{source-lb}, we
immediately obtain:
\begin{remark}\label{source-npir}
For the NPIR models, suppose that Assumption \ref{npir-a1} holds,
and Assumption \ref{h-smooth} is satisfied with $h\in {\cal
H}_{source}$ and $u_k = e_k$, $\nu_k^{-r} = \psi (\lambda_k)$ for
all $k\ge 1$. Let $\phi(B^{-2})=K^\ast K$. Then the estimator $\hat
h_n$ defined in (\ref{npir-estimate}) with $u_k = e_k$ reaches the
minimax rate uniformly over $h\in {\cal H}_{source}$:
\[
\sup_{h\in {\cal H}_{source}}\E_h[\norm{\hat h_n-h}_X^2]\le
C\delta_n,\quad \text{with}~\delta_n ~\text{and}~m~\text{given in
Remark \ref{source-lb}}.
\]
\end{remark}
In the literature on ill-posed inverse problems with known operator
$K$, there are available many other estimation procedures (like
Tikhonov's method) that employ source conditions; some of which lead
to rate-optimal estimators only for mildly ill-posed case. See,
e.g., Bissantz, Hohage, Munk and Ruymgaart (2007) and Florens,
Johannes and Van Bellegem (2007) for recent results.

\subsection{Relaxation of the exact link condition}

For the minimax risk lower bound we impose Assumption
\ref{K-smooth}, and for the upper bound we use Assumption
\ref{npir-a2}. Together, these two assumptions require that the
operator $K$ satisfies
\[
c\norm{[\phi(B^{-2})]^{1/2}h}_X\le \norm{K h}_W\le M
\norm{[\phi(B^{-2})]^{1/2}h}_X \text{ for all $h\in L_X^2$},
\]
which is equivalent to
\begin{equation}
\label{exact-range} \ran( [\phi(B^{-2})]^{1/2}) =\ran (\abs{K} ).
\end{equation}%
This is a standard condition imposed even in books and papers on
ill-posed inverse problems with deterministic errors; see, e.g.,
Engl, Hanke and Neubauer (1996), Nair, Pereverzev and Tautenhahn
(2005) and the references therein. This condition usually holds when
$K$ acts exactly along certain function classes; see Section
\ref{sec-app} for such an example. Moreover, this exact range
condition \ref{exact-range} is automatically satisfied under the
source condition with $K^\ast K = \phi(B^{-2})$. However, Assumption
\ref{npir-a2} may fail more generally. Luckily, this assumption is
not strictly necessary.

Let us indicate one possibility how Assumption \ref{npir-a2} can be
relaxed to requiring
\[
\ran [\phi(B^{-2})]^{1/2}\subset \ran \abs{K}+L, ~\text{for some
finite-dimensional linear space}~ L.
\]
To keep it simple, we consider the case that the subspace $L$ is
spanned by one eigenfunction $u_\ell$ of $B$ with $u_\ell\notin \ran
\abs{K}$ and $1\le \ell \le m$. Then the simple estimator $\hat h_n$
using $\hat\eta_\ell$ given in \eqref{npir-estimate} is no longer
well defined, but we can consider for some $v\in L^2_W$ the
estimator
\[ \tilde\eta_\ell:=\frac1n\sum_{i=1}^nY_i v(W_i),\quad \tilde h_n : =\sum_{k=1, k\neq \ell}^m\hat\eta_k u_k + \tilde\eta_\ell u_\ell.\]
Following the bias variance decomposition in the proof of
Proposition \ref{prop-npir}, we obtain
\begin{align}
\E[(\tilde\eta_\ell-\scapro{h}{u_\ell}_X)^2]&\le(\scapro{Kh}{v}_W-\scapro{h}{u_\ell}_X)^2+2n^{-1}(S^2+\sigma_1^2)\norm{v}_W^2
\nonumber\\
& \asymp \scapro{h}{K^\ast
v-u_\ell}_X^2+n^{-1}\norm{v}_W^2.\label{eq-coeffmseh}
\end{align}
The definition of $H^r_R$ implies with some uniform constant $C>0$
\begin{equation}\label{eq-coeffmsemax}
\sup_{h\in {H^r_R}}\E[(\tilde\eta_\ell-\scapro{h}{u_\ell}_X)^2]\le
C(R^2\norm{B^{-r}(K^\ast v-u_\ell)}_X^2+n^{-1}\norm{v}_W^2).
\end{equation}%
From inequality (\ref{eq-coeffmsemax}), it is easy to derive that
this error in estimating the coefficient $\scapro{h}{u_\ell}_X$ is
minimized by
\[ v=(KK^\ast+n^{-1}R^{-2}B^{2r})^{-1}Ku_\ell,\]
which is always well-defined. Consequently, in terms of minimax
optimal rate over the class of functions $H^r_R$, the rate in
Proposition \ref{prop-npir} does not deteriorate if we use
$\tilde\eta_\ell$ instead of $\hat\eta_\ell$ and its error bound
\[n^{-1}\norm{(KK^\ast+n^{-1}R^{-2}B^{2r})^{-1}Ku_\ell}_W^2\]
is not larger than the minimax optimal rate. See Section
\ref{sec-app} for a concrete example.

\section{An upper bound for the NPIV model}

We now provide an upper bound for the NPIV model. For the NPIV model
additional considerations due to the unknown conditional expectation
operator are necessary. It is, of course, more complex to construct
an estimator that is rate-optimal for the NPIV model than for the
NPIR model, which is why the approaches in the literature are more
diverse and require different additional assumptions. Here, we
restrict ourselves to presenting a simple estimator to illustrate
that it is possible to construct a rate-optimal estimator for the
NPIV model in both mildly ill-posed and severely ill-posed cases
based on the SMD estimator of Newey and Powell (2003), Ai and Chen
(2003) and Blundell, Chen and Kristensen (2007). First, for each
integer $J \ge 1$, we denote by $\spann\{p_1,...,p_J\}$ a
$J$-dimensional linear subspace of $L^2_W$ that becomes dense in
$L^2_W$ as $J\to\infty$. Let $P^{J_n} (w)=(p_1 (w),...,p_{J_n}
(w))^{\prime}$ and $\mathbf{P}
=(P^{J_{n}}(W_{1}),...,P^{J_{n}}(W_{n}))^{\prime }$. We compute a
sieve least squares estimator of $E[Y-h(X)|W=\cdot ]$ as
\[
\widehat{\E}[Y-h(X)|W=\cdot]=\sum_{t=1}^{n}\{Y_{t}-h(X_{t})\}P^{J_{n}}(W_{t})^{\prime }(\mathbf{P}%
^{\prime }\mathbf{P})^{-1}P^{J_{n}}(\cdot).
\]
For each integer $m\ge 1$, we denote by $\mathcal{H}_m
:=\spann\{\psi_1,...,\psi_m\}$ an $m$-dimensional linear subspace of
$L^2_X$ that becomes dense in $L^2_X$ as $m\to\infty$. Then we
compute the SMD estimator of the true structure function $h_0$ as
\begin{equation}
\label{npiv-estimate}
\hat{h}_n =\argmin_{h\in \mathcal{H}_{m(n)}\cap \cal H (r,R)}\frac{1}{n}\sum_{i=1}^{n}%
\Big\{\widehat{\E}[Y-h(X)|W=W_i ]\Big\}^{2}.
\end{equation}%
Depending on the prior information about ${\cal H}(r,R)$, sometimes
one may compute $\hat{h}_n$ in closed form. For example, if ${\cal
H}(r,R) = H^r_R$ and the density of $X$ is bounded below and above
by positive  constants, then
\begin{align}
\hat h_n (x)&=\sum_{k=1}^m\hat\pi_k \psi_k(x)
= \psi^{m}(x)^{\prime }\widehat{\Pi },\\
\widehat{\Pi }&=\left( \Psi ^{\prime }\mathbf{P}(\mathbf{P}^{\prime }\mathbf{P})^{-1}%
\mathbf{P}^{\prime }\Psi +\hat{\lambda} \mathbf{C}\right) ^{-1}\Psi ^{\prime }\mathbf{P}(%
\mathbf{P}^{\prime }\mathbf{P})^{-1}\mathbf{P}^{\prime }\mathcal{Y},
\end{align}%
with $\Psi =(\psi ^{m}(X_{1}),...,\psi ^{m}(X_{n}))^{\prime }$,
$\mathcal{Y} =(Y_{1},...,Y_{n})^{\prime }$, the penalization matrix
$\mathbf{C} = \int \{[B^r \psi ^{m}(x)] [B^r \psi ^{m}(x)]^{\prime
}\}dx$ and $\hat{\lambda }$ satisfies $\widehat{\Pi}^{\prime
}\mathbf{C}\widehat{\Pi } = R^2$.

In addition to the assumptions on the NPIR models, we impose the
following:
\begin{assumption}
\label{npiv-a1} The basis $\{\psi_k\}$ is a Riesz basis associated
with the operator $B$, that is, $\sum_{k=1}^{\infty}
\scapro{h}{\psi_k}_X^2 \asymp \sum_{k=1}^{\infty}
\scapro{h}{u_k}_X^2$ for all $h\in L_X^2$.
\end{assumption}
Assumption \ref{npiv-a1} allows for the use of a Riesz basis
$\{\psi_k\}$ instead of the ideal orthonormal basis $\{u_k\}$ to
approximate the unknown structure function $h\in {\cal H} (r,R)$
with the same order of the approximation errors. Of course in
applications, we need some information about the tail behavior of
the density of $X$ before we can construct such a basis. For
example, if we know that the density of $X$ is bounded above and
below by finite positive constants , then we could use wavelets as
the $\{\psi_k\}$.
\begin{assumption}\label{npiv-a2}\mbox{}\\
(1) $\E[Y-\Pi_m (h(X))|W=\cdot ]$ belongs to $\Lambda
_{c}^{r_{K}}(\cal W)$ (H\"older ball of regularity $r_K$)
for any $\Pi_m (h) \in \mathcal{H}_{m}$;\\
(2) (i) the smallest and the largest eigenvalues of $%
\E\{P^{J_n}(W)P^{J_n}(W)^{\prime }\}$ are bounded and bounded away
from zero for each $J_{n}$; (ii) $P^{J_n}(W)$ is a tensor product of
either a cosine series or a B-spline basis of order $\gamma _{b}$ or
a wavelet basis of order $\gamma_{b}$, with $\gamma
_{b}>r_{K}> d/2$;\\
(3) the density of $W$\ is continuous and bounded away from zero
over its support $\cal W$, which is a
compact connected subset in $\R^d$ with Lipschitz continuous boundaries and non-empty interior;\\
(4) (i) $J_n\rightarrow \infty $ and $J^2_n /n \rightarrow 0 $; (ii)
$\lim_{n}\frac{J_n}{m(n)}=c \in (1, \infty)$ and $J_n > m(n)$.
\end{assumption}
Assumption \ref{npiv-a2} implies that the sieve least square
estimate $\widehat{\E}[h(X)|W=\cdot]$ of $\E[h(X)|W=\cdot]$ performs
well; see e.g., Blundell, Chen and Kristensen (2007) for details.

\begin{theorem}\label{thm-npiv}
For the NPIV models, suppose that Assumptions \ref{h-smooth},
\ref{K-smooth}, \ref{npir-a1}, \ref{npir-a2}, \ref{npiv-a1} and
\ref{npiv-a2} hold. Then the estimator $\hat h_n$ defined in
(\ref{npiv-estimate}) satisfies
\[
\norm{\hat h_n-h}_X^2 \le C \max\Big\{\nu_{m+1}^{-2r},~ \frac{m}{n}
[\phi(\nu_{m}^{-2})]^{-1}\Big\}
\]
uniformly over $h\in {\cal H} (r,R)$ except on an event whose
probability tends to zero as $n\uparrow \infty$.
If $m=m(n)$ is such
that $n^{-1}\sum_{k=1}^m\nu_k^{2r}[\phi(\nu_k^{-2})]^{-1}\asymp 1$,
then this estimator $\hat h_n$ is rate-optimal in the minimax sense:
there is a finite constant $C>0$ such that
\[
\norm{\hat h_n-h}_X^2\le C\nu_{m+1}^{-2r} \asymp \frac{m}{n}
[\phi(\nu_{m}^{-2})]^{-1} \asymp \delta_n,\quad \text{with}~\delta_n
~\text{given in Theorem \ref{thm-lb}},
\]
uniformly over $h\in {\cal H} (r,R)$ except on an event whose
probability tends to zero as $n\uparrow \infty$.

(1) Mildly ill-posed case: Let $\phi(t)=t^a$ and $\nu_k\asymp
k^\epsilon$ for some $a, \epsilon >0$. If $m\asymp n^{1/(2r\epsilon
+2a\epsilon+1)}$, then $\delta_n \asymp
n^{-2r/(2r+2a+\epsilon^{-1})}$.

(2) Severely ill-posed case: Let $\phi(t)=\exp(-t^{-a/2})$,
$\nu_k\asymp k^\epsilon$ for some $a, \epsilon>0$. If
$m=c\log(n)^{1/a\epsilon}$ with a sufficiently small $c>0$, then
$\delta_n \asymp (\log n)^{-2r/a}$.
\end{theorem}

This minimax rate theorem appears to be new in the literature, and
can be proved by slightly modifying the proof of Blundell, Chen and
Kristensen (2007) for their theorem 2. Hall and Horowitz (2005)
obtained minimax optimal rate $\sup_{h\in {\cal
H}_{HH}}\E_h[\norm{\hat h_n-h}_X^2]\le C\delta_n$ for their
estimators in the mildly ill-posed case for the class of functions
${\cal H}_{HH}$ defined in (\ref{HH}). Hoffmann and Rei\ss\ (2007)
propose a wavelet estimator in the case of an unknown operator $K$
that is elliptic with ill-posedness degree $a$. They assume there
exists an estimator of $K$ with specified rate, and their class of
functions ${\cal H} (r,R)$ is a Besov ball that could be bigger than
the function class defined in our Assumption \ref{h-smooth}, but
they do not consider severely ill-posed case.

\section{More on regularity conditions}\label{sec-app}

In this section, we use examples to discuss the pros and cons of the
approach of imposing two basic regularity conditions (the
approximation and the link conditions) versus the other approach of
using the general source condition. To simplify the discussion, here
we assume the operator $K$ is known. In the first class of examples,
the operator $K$ has very smooth eigenfunction basis (in the sense
that its eigenfunctions are many times differentiable), while in the
second class of examples, the operator $K$ has eigenfunctions that
are not differentiable.

\subsection{Examples of $K$ having infinitely times differentiable eigenfunctions}

Suppose that the $\{W_i\}_{i=1}^{n}$ are uniformly distributed on
$[0,1]$ and $K$ is a circular convolution operator on $L^2([0,1])$:
$Kh(w)=\int_0^1 k(x-w)h(x)dx$ with a $1$-periodic function $k$ that
satisfies $k(-x)=k(x)$ and has Fourier coefficients $\abs{{\cal
F}k(m)}=\abs{\int_0^1 k(x)\cos(mx)\,dx}\asymp (1+\abs{m})^{-a}$.
Then $K$ is a positive-definite self-adjoint operator which is
diagonalized by the Fourier basis.
%In statistical terms this
%prescription means that $X_i=W_i+V_i \mod 1$ (explain mod?\MR) with
%some independent random variable $V_i\in [0,1]$ of density $k$. Note
%that then also the $X_i$ are uniform on $[0,1]$.
\begin{description}
\item[Source condition:] In this canonical case the exact link between $K$ and $B$ is easily
established with $B=K^{-1}$, $\phi(t)=t$ hence
$\norm{Kh}_W=\norm{[\phi(B^{-2})]^{1/2}h}_X$ for all $h\in
L^2([0,1])$. The smoothness of the unknown function $h$ is also
described using $B=K^{-1}$; hence the Hilbert scale space $H^r$
(generated by $B$) is equal to the classical periodic Sobolev space
$H^{ra}_{per}$ of smoothness (or regularity) $ra$. Applying Remark
\ref{source-npir}, we obtain minimax optimal rate for this scale of
periodic Sobolev spaces.

\item[Approximation + link conditions:] Suppose $\{u_k\}_{k\ge 1}$ is an orthonormal basis of
$L^2([0,1])$ such that $\norm{Kg}_{L^2([0,1])}^2\asymp
\sum_{k=1}^\infty k^{-2a}\scapro{g}{u_k}^2$. A typical example is
given by sufficiently regular periodized wavelet bases (see Cohen,
Daubechies and Vial (1993)). Then we can define
\[Bg:=\sum_{k\ge 1}k\scapro{g}{u_k}u_k,\]
and the Hilbert scale spaces $H^r$ can be interpreted as
approximation spaces for the basis $(u_k)$. In the convolution
example we obtain $\norm{Kg}_W\asymp \norm{B^{-a}g}_X$.
Consequently, the exact link conditions (assumptions \ref{K-smooth}
and \ref{npir-a2}) between $K$ and $B$ hold with $\phi(t)=t^a$.
Applying Proposition \ref{prop-npir}, we obtain minimax optimal rate
for the Hilbert scale space $H^r$ generated by $B$.
\end{description}

The Hilbert scale of approximation spaces generated by $B$ does not
necessarily coincide with the Hilbert scale generated by $K$. The
most pronounced example is the case $a<1/2$, where all non-periodic
wavelets on an interval still satisfy
$\norm{Kg}_{L^2([0,1])}^2\asymp \sum_{k=1}^\infty
k^{-2a}\scapro{g}{u_k}^2$ (see Cohen, Daubechies and Vial (1993)).
Hence, the approximation spaces for unknown true structure function
need not exhibit any boundary condition. This means that a smooth,
but non-periodic function on $[0,1]$ will have high regularity $r$
in terms of the approximation space, while it is an element in
periodic Sobolev spaces up to regularity $1/2$ only. If we have in
mind that our true function $h$ is smooth, but not periodic, we
should therefore rather choose the approximation space approach. On
the other hand, wavelets work well just to some maximal regularity
and they will therefore reconstruct very smooth and periodic
functions not as well as the Fourier basis.

If $K$ is more ill-posed, that is $a\ge 1/2$, we can adopt the ideas
explained in Subsection 4.2. We remain in the approximation space
framework and use non-periodic compactly supported wavelets as basis
functions $\{u_k\}$. Only the wavelets $\psi_\lambda$ with support
at the boundary are not in the periodic Sobolev spaces $H^s_{per}$,
$s\ge 1/2$. Using some (statistical) kernel function $L_h:
[-b,b]\to\R$ of bandwidth $b$, we can consider the periodically
smoothed version
\[ \tilde \psi_\lambda(x):=\int_{-b}^b
\psi_\lambda(\{x-y\})L_b(y)\,dy,\quad x\in [0,1],
\]
where $\{z\}=z-\floor{z}\in [0,1)$ denotes the fractional part of
$z\in\R$. If $L$ and $\psi_\lambda$ are sufficiently often
differentiable, then $\tilde\psi_\lambda$ lies in the range
$H^a_{per}$ of $K$. Using $v=K^{-1}\tilde\psi_\lambda$ in equation
\eqref{eq-coeffmseh}, standard kernel estimates ($h\in H^r_R$
implies $h\in H^s_{per}$ for all $s\le r$ and $s<1/2$) show that for
all $h\in H^r_R$ (with adapted notation)
\[\E[(\tilde\eta_\lambda-\scapro{h}{\psi_\lambda})^2]\le
C_1(\scapro{h}{\tilde\psi_\lambda-\psi_\lambda}^2+n^{-1}\norm{K^{-1}\tilde\psi_\lambda}^2)
\le C_2(b^{2s}+n^{-1}b^{-2a}).
\]
Optimizing over $b$, we infer that $\scapro{h}{\psi_\lambda}$ can be
estimated at rate $n^{-s/(s+a)}$, which for $r\ge 1/2$ is nearly
$n^{-1/(2a+1)}$. Since in a wavelet approximation space of dimension
$2^J$ only of the order $J$ wavelets lie at the boundary, the rate
in estimating $h$ will be $n^{-s/(s+a)}\log(n)+n^{-2r/(2r+2a+1)}$,
which for $r\ge \frac12+\frac{1}{4a}$ is roughly $n^{-1/(2a+1)}$. If
we had taken a method based on the source condition approach (like
projection on eigenfunctions of $K$, or Tikhonov methods) the best
achievable rate would have been roughly $n^{-1/(2a+2)}$.

\subsection{Examples of $K$ having non-differentiable eigenfunctions}

Depending on applications, it is perfectly conceivable that the
eigenfunctions of $K$ are rough while the basis functions $u_k$ of
$B$ are smooth (or differentiable). For example, we can use the Haar
basis $\psi_{jk}(x)=\psi(2^jx-k)$ on $L^2([0,1])$ ($\psi(x)={\bf
1}_{[0,1/2]}-{\bf 1}_{[1/2,1]}$, $j\in\N_0$, $k=0,\ldots,2^j-1$, and
$\psi_{-1,0}={\bf 1}_{[0,1]}$) and define -- somewhat artificially
-- in this Haar basis
\[ K\psi_{jk}:=2^{\alpha j}\psi_{jk}.\]
Then $K$ is self-adjoint with eigenfunctions $\psi_{jk}$, which are
step functions. For $\alpha r<1/2$, the Hilbert scale $H^r$ of $K$
(or of the Harr basis) will be a Sobolev space, whereas for any
$\alpha r\ge 1/2$ this Hilbert scale $H^r$ will not be described in
terms of traditional smoothness. Note that this $H^r$ will always
contain piecewise constant jump functions. Nevertheless, the larger
$r$ the less complex is the function class ${\cal H}(r,R)$, that is
the smaller the approximation error rate. As for the convolution
operator we could instead define the function class ${\cal H}(r,R)$
in terms of a basis $\{u_k\}$ associated to $B$ which is smoother
and satisfies at the same time the link conditions of Assumptions
\ref{K-smooth} and \ref{npir-a2}.

In conclusion, we see that rate-optimal methods may behave poorly if
the function of interest, the structure function $h$, is not regular
in the setting for which the method is designed. An important part
of the specification of rate optimality is therefore always the
associated function class.

\section{Perspectives}

In this paper, we clarify the relations between the existing sets of
regularity conditions for convergence rates of NPIV regression
models. We establish minimax risk lower bounds in mean squared error
loss for the NPIV models under two basic regularity conditions that
allow for both mildly ill-posed and severely ill-posed cases. We
also show that the simple SMD estimator achieves the minimax risk
lower bound, hence is rate-optimal for both mildly ill-posed and
severely ill-posed cases.

Many of the ideas in this paper can be easily adapted to treat other
kinds of ill-posed inverse problems in econometrics. For instance,
when the problem is mildly ill-posed, Horowitz and Lee (2007) show
that their kernel based Tikhonov regularized estimator of
nonparametric quantile instrumental variables (IV) regression
reaches the minimax rate under conditions very similar to those
imposed in Hall and Horowitz (2005) for NPIV regression. Similarly,
one could show that the penalized SMD estimator proposed in Chen and
Pouzo (2007) for nonlinear and possibly nonsmooth nonparametric
conditional moment models is also rate-optimal, as their estimator
achieves the minimax risk lower bounds established in our paper for
the NPIV regression model.

Once this is established, the intriguing open problem remains how to
choose the regularization parameters adaptively from the data, not
knowing the true regularity, and even to select among the different
proposed procedures (e.g. generated by different operators $B$) in a
data-driven way.

%\thebibliography{99}

%\harvarditem{Bissantz {\em et al.}}{2007}{BHMR} {\sc Bissantz, N.,
%Hohage, T., Munk, A., Ruymgaart, F.} (2007). Convergence rates of
%general regularization methods for statistical inverse problems and
%applications. {\em SIAM J. Numerical Analysis}, to appear.

%\harvarditem{Cohen, Daubechies and Vial}{1993}{CDV} {\sc Cohen, A.,
%Daubechies, I. and Vial P.} (1993). Wavelets on the interval and
%fast wavelet transforms. {\em Appl. Comp. Harm. Anal.} {\bf 1}(1),
%54--81.

%\harvarditem{Korostelev and Tsybakov}{1993}{KoTs} {\sc Korostelev,
%A. and Tsybakov, A.} (1993). {\em Minimax theory of image
%reconstruction.} Lecture Notes in Statistics 82, Springer, New York.

\bigskip

\thebibliography{99}

%\begin{description}
\item \textsc{Ai, C. and X. Chen} (2003): \textquotedblleft Efficient
Estimation of Models with Conditional Moment Restrictions Containing
Unknown Functions,\textquotedblright\ \emph{Econometrica}, 71,
1795-1844.

\item \textsc{Bissantz, N., T. Hohage, A. Munk and F. Ruymgaart} (2007). Convergence rates of
general regularization methods for statistical inverse problems and
applications. {\em SIAM J. Numerical Analysis}, to appear.

\item \textsc{Blundell, R., X. Chen and D. Kristensen} (2007):
\textquotedblleft Semi-nonparametric IV Estimation of
Shape-Invariant Engel Curves,\textquotedblright\ forthcoming in
\textit{Econometrica}.

\item \textsc{Blundell, R. and J.L. Powell} (2003): \textquotedblleft
Endogeneity in Semiparametric and Nonparametric Regression
Models,\textquotedblright\ In: \textit{Advances in Economics and
Econometrics: Theory and Applications}, Dewatripont, M., Hansen,
L.P., and Turnovsky, S.J., eds, vol. 2, pp. 312-357. Cambridge, UK:
Cambridge University Press.

\item \textsc{Carrasco, M., J.-P. Florens and E. Renault} (2007):
\textquotedblleft Linear Inverse Problems in Structural Econometrics
Estimation Based on Spectral Decomposition and
Regularization\textquotedblright ,\ in J.J. Heckman and E.E. Leamer
(eds.), \textit{The Handbook of Econometrics}, vol. 6.\
North-Holland, Amsterdam, forthcoming.

\item \textsc{Chen, X. and D. Pouzo} (2007): \textquotedblleft On
Estimation of Semi/nonparametric conditional moment models with
possibly nonsmooth moments,\textquotedblright\ manuscript, Yale
University and New York University, Dept. of Economics.

\item \textsc{Chernozhukov, V. and C. Hansen} (2005):
\textquotedblleft An IV Model of Quantile Treatment Effects,
\textit{Econometrica}, 73, 245-261.

\item \textsc{Chernozhukov, V., G. Imbens, and W. Newey} (2007):
\textquotedblleft Instrumental Variable Estimation of Nonseparable
Models,\textquotedblright\ forthcoming in \textit{Journal of
Econometrics}.

\item \textsc{Cohen, A., I. Daubechies, and P. Vial} (1993). Wavelets on the interval and
fast wavelet transforms. {\em Appl. Comp. Harm. Anal.} {\bf 1}(1),
54--81.

\item \textsc{Cohen, A., M. Hoffmann, and M. Rei\ss }\ (2004):
\textquotedblleft Adaptive wavelet Galerkin methods for linear
inverse problems,\textquotedblright\ \textit{SIAM J. Numer. Anal.}
42, 1479-1501.

\item \textsc{Darolles, S., J.-P. Florens and E. Renault} (2002):
\textquotedblleft Nonparametric Instrumental
Regression,\textquotedblright\ mimeo, GREMAQ, University of
Toulouse.

\item \textsc{Donoho, D., R. Liu and B. MacGibbon} (1990): \textquotedblleft Minimax Risk Over Hyperrectangles, and
Implications,\textquotedblright\ \textit{Annals of Statistics.} 18,
1416-1437.

\item \textsc{Edmunds, D. and W. Evans} (1987): \textit{Spectral Theory and
Differential Operators}, Oxford University Press: Oxford.

\item \textsc{Efromovich, S. and V. Koltchinskii} (2001): \textquotedblleft
On Inverse Problems with Unknown Operators,\textquotedblright\
\textit{IEEE Trans. on Information Theory}, 47, 2876-2893.

\item \textsc{Engl, H., M. Hanke and A. Neubauer} (1996): \textit{%
Regularization of Inverse Problems}, Kluwer Academic Publishers:
London.

\item \textsc{Florens, J.-P. }(2003): \textquotedblleft Inverse Problems and
Structural Econometrics: the Example of Instrumental
Variables,\textquotedblright\ In: \textit{Advances in Economics and
Econometrics: Theory and Applications}, Dewatripont, M., Hansen,
L.P., and Turnovsky, S.J., eds, vol. 2, pp. 284-311. Cambridge, UK:
Cambridge University Press.

\item \textsc{Florens, J.-P., J. Johannes and S. Van Bellegem}
(2007): \textquotedblleft Identification and Estimation by
Penalization in Nonparametric Instrumental
Regression,\textquotedblright\ mimeo, GREMAQ, University of
Toulouse.

\item \textsc{Gagliardini, P. and O. Scaillet} (2006):\textquotedblleft Tikhonov Regularisation for Functional
Minimum Distance Estimation,\textquotedblright\ mimeo
%
%\item {\sc Goldenshluger, A. and S. V. Pereverzev} (2003): On adaptive inverse
%estimation of linear functionals in Hilbert scales. {\sl Bernoulli}
%9(5), 783--807.

\item \textsc{Hall, P. and J. Horowitz }(2005): \textquotedblleft
Nonparametric Methods for Inference in the Presence of Instrumental
Variables,\textquotedblright\ \emph{Annals of Statistics}, 33,
2904-2929.

\item \textsc{Hoffmann, M. and M. Rei\ss }\ (2007):
\textquotedblleft Nonlinear estimation for linear inverse problems
with error in the operator,\textquotedblright\ \textit{Annals of
Statistics}, to appear.

\item \textsc{Horowitz, J. and S. Lee }(2007): \textquotedblleft
Nonparametric Instrumental Variables Estimation of a Quantile
Regression Model,\textquotedblright\ forthcoming in
\textit{Econometrica}.

\item \textsc{Korostelev, A. and A. Tsybakov }(1993): \textit{Minimax Theory
of Image Reconstruction}, Lecture Notes in Statistics, v.82.
Springer: New York.

\item \textsc{Le Cam, L. and G. Lo Yang} (2000): \textit{Asymptotics in statistics. Some basic concepts.} 2nd
ed., Springer, New York.

\item \textsc{Nair, M., S. Pereverzev and U. Tautenhahn }(2005): \textquotedblleft
Regularization in Hilbert scales under general smoothing
conditions,\textquotedblright\ \emph{Inverse Problems}, 21,
1851-1869.

\item \textsc{Newey, W.K. and J. Powell }(2003): \textquotedblleft
Instrumental Variables Estimation for Nonparametric
Models,\textquotedblright\ \emph{Econometrica}, 71, 1565-1578.

\item \textsc{Yang, Y. and A. Barron} (1999): Information-theoretic
determination of minimax rates of convergence, {\em Ann. Stat.}
27(5), 1564--1599.
%\end{description}

\bigskip

\begin{center}
{\Large \textbf{Appendix: Proofs}}

\medskip
\end{center}

\begin{proof}[Proof of Lemma \ref{LemReduction}]

Let $\hat h_n=\hat h_n(\{(X_i,Y_i,W_i)\}_{i=1}^{n} )$ be an
estimator for the NPIV model. Knowing the operator $K$ amounts to
knowing the conditional law of $X_i$ given $W_i$. Let us call the
observations in the NPIR model $\{(Y_i',W_i')\}_{i=1}^{n}$ for some
$({\cal L}_W,\sigma(\cdot),h)\in{\cal C_0}$. We then generate
artificially i.i.d. observations $X_i'$ according to the conditional
law ${\cal L}_{X|W=w}$ with ${w=W_i'}$. Then the observations
$\{(X_i',Y_i',W_i')\}_{i=1}^{n}$ follow the law of some $({\cal
L}_{UWX},h)\in{\cal C}$ because $Y_i'=h(X_i')+U_i'$ holds with
$U_i'=(Kh)(W_i')-h(X_i')+V_i'$ satisfying $\E[U_i'\,|\,W_i']=0$ and
${\cal L}_{V'|W'=w}=N(0,\sigma^2(w))$. Consequently, the
(randomized) estimator $\tilde h_n(\{(Y_i',W_i')\}_{i=1}^{n}):=\hat
h_n(\{(X_i',Y_i',W_i')\}_{i=1}^{n})$ has the same risk under $({\cal
L}_W,\sigma(\cdot),h)\in{\cal C_0}$ as $\hat h_n$ has under $({\cal
L}_{U'W'X'},h)\in{\cal C}$, and is thus not larger than the maximal
risk over $\cal C$.
\end{proof}

\begin{proof}[Proof of Theorem \ref{thm-lb}]
We consider for $\theta=(\theta_k)$ with $\theta_k\in\{-1,+1\}$ and
a sequence $(\gamma_k)$, to be specified below, the following
functions in $L_X^2$:
\[ h_\theta:=\sum_{k=1}^m \theta_k\gamma_k u_k.\]
The property $h_\theta\in H_R^r$ yields the following constraint on
$m$ and $(\gamma_k)$:
\[ \norm{h_\theta}_r^2=\sum_{k=1}^m\nu_k^{2r}\gamma_k^2\le R^2.
\]
For $\ell=1,\ldots,m$ and each $\theta$ introduce $\theta^{(\ell)}$
by $\theta^{(\ell)}_k=\theta_k$ for $k\not=\ell$ and
$\theta^{(\ell)}_\ell=-\theta_\ell$. Then because of the Gaussianity
of the $V_i$ given $W_i$ the log-likelihood of
$\PP_{\theta^{(\ell)}}$ w.r.t. $\PP_\theta$ is
\[ \log\Big(\frac{d\PP_{\theta^{(\ell)}}}{d\PP_\theta}\Big)=\sum_{i=1}^n\pm
\frac{2\gamma_\ell (K u_\ell)(W_i)}{\sigma^2(W_i)}V_i
-\frac12\sum_{i=1}^n \Big( \frac{2\gamma_\ell (K
u_\ell)(W_i)}{\sigma(W_i)}\Big)^2.
\]
Its expectation satisfies
\begin{align*}
\E_\theta\Big[\log\Big(\frac{d\PP_{\theta^{(\ell)}}}{d\PP_\theta}\Big)\Big]&=
-2\gamma_\ell^2n\norm{(Ku_\ell)\sigma^{-1}}_W^2\\
&\ge
-2M\sigma_0^{-2}\gamma_\ell^2n\norm{[\phi(B^{-2})]^{1/2}u_\ell}_X^2\\
&=-2M\sigma_0^{-2}\gamma_\ell^2n\phi(\nu_\ell^{-2})=:\mu_\ell.
\end{align*}
In terms of the Kullback-Leibler divergence this means
$\KL(\PP_{\theta^{(\ell)}},\PP_\theta)\le -\mu_\ell$. More
explicitly, we obtain by Markov's inequality
\[\PP_\theta\Big(-\frac12\sum_{i=1}^n
\Big( \frac{2\gamma_\ell (K u_\ell)(W_i)}{\sigma(W_i)}\Big)^2\le
-2\mu_\ell\Big)\le \frac{-\mu_\ell}{-2\mu_\ell}=\frac12.
\]
Using the symmetry of the distribution of $V_i$ given $W_i$, we
infer by conditioning on $(W_i)_{1\le i\le n}$
\[
\PP_\theta\Big(\frac{d\PP_{\theta^{(\ell)}}}{d\PP_\theta}\ge
\exp(2\mu_\ell)\Big)
=\E_\theta\Big[\PP_\theta\Big(\log\Big(\frac{d\PP_{\theta^{(\ell)}}}{d\PP_\theta}\Big)
\ge 2\mu_\ell\,\Big|\,(W_i)_{1\le i\le n} \Big)\Big]\ge \frac12.
\]
We calculate for each estimator $\hat h_n$:
\begin{align*}
&\sup_{h\in H_R^r}\E_h[\norm{\hat h_n-h}_X^2]\\
&\ge
\sup_{\theta\in\{-1,+1\}^m}\E_{\theta}[\norm{\hat h_n-h_\theta}_X^2]\\
&\ge
2^{-m}\sum_{\theta\in\{-1,+1\}^m}\sum_{k=1}^m\E_{\theta}[\scapro{\hat
h_{n}-h_{\theta}}{u_k}_X^2]\\
&=
\sum_{k=1}^m2^{-m}\sum_{\theta\in\{-1,+1\}^m}\frac12\Big(\E_{\theta}[\scapro{\hat
h_{n}-h_{\theta}}{u_k}_X^2]+\E_{{\theta^{(k)}}}[\scapro{\hat
h_{n}-h_{\theta^{(k)}}}{u_k}_X^2]\Big)\\
&=
\sum_{k=1}^m2^{-m}\sum_{\theta\in\{-1,+1\}^m}\frac12\E_{\theta}\Big[\scapro{\hat
h_{n}-h_{\theta}}{u_k}_X^2+\scapro{\hat
h_{n}-h_{\theta^{(k)}}}{u_k}_X^2\frac{d\PP_{\theta^{(k)}}}{d\PP_\theta}\Big]\\
&\ge
\sum_{k=1}^m2^{-m}\sum_{\theta\in\{-1,+1\}^m}\frac{\exp(2\mu_k)}{2}\E_{\theta}\Big[\Big(\scapro{\hat
h_{n}-h_{\theta}}{u_k}_X^2+\scapro{\hat
h_{n}-h_{\theta^{(k)}}}{u_k}_X^2\Big)\times\\
&\qquad \times{\bf
1}\big\{\tfrac{d\PP_{\theta^{(k)}}}{d\PP_\theta}\ge \exp(2\mu_k)\big\}\Big]\\
&\ge
\sum_{k=1}^m2^{-m}\sum_{\theta\in\{-1,+1\}^m}\frac{\exp(2\mu_k)}{8}\scapro{h_{\theta}-h_{\theta^{(k)}}}{u_k}_X^2
\PP_\theta\Big(\frac{d\PP_{\theta^{(k)}}}{d\PP_\theta}\ge \exp(2\mu_k)\Big)\\
&\ge  \sum_{k=1}^m \frac{\exp(2\mu_k)}{4} \gamma_k^2.
\end{align*}
We choose $\gamma_k=\sigma_0n^{-1/2}[\phi(\nu_k^{-2})]^{-1/2}$ such
that $\mu_k=-2M$ and then pick the largest $m\ge 1$ such that
$\sum_{k=1}^m\nu_k^{2r}\gamma_k^2\le R^2$.

This gives the lower bound
\[
\inf_{\hat h_n}\sup_{h\in{\cal H}(r,R)}\E_h[\norm{\hat
h_n-h}_X^2]\ge \inf_{\hat h_n}\sup_{h\in H_R^r}\E_h[\norm{\hat
h_n-h}_X^2]\ge \tfrac{\sigma_0^2}{4\exp(4M)} n^{-1}\sum_{k=1}^m
[\phi(\nu_k^{-2})]^{-1}
\]
where $m$ is largest possible with
$\sum_{k=1}^m\nu_k^{2r}\gamma_k^2\le R^2$, i.e.
\[\sigma_0^2n^{-1}\sum_{k=1}^m\nu_k^{2r}[\phi(\nu_k^{-2})]^{-1}\le R^2.
\]
(1) (mildly ill-posed case): When $\phi(t)=t^a$ and $\nu_k\asymp
k^\epsilon$ for some $a, \epsilon>0$, we have asymptotically as
$n\to\infty$:
\[
n^{-1}\sum_{k=1}^m\nu_k^{2r}[\phi(\nu_k^{-2})]^{-1}=n^{-1}\sum_{k=1}^m
k^{2\epsilon r+2\epsilon a}\asymp n^{-1}m^{2\epsilon r+2\epsilon
a+1}.
\]
Hence, choosing $m\asymp n^{1/(2\epsilon r+2\epsilon a+1)}$ we
obtain the asymptotic lower bound
\[\delta_n \asymp n^{-1}\sum_{k=1}^m
[\phi(\nu_k^{-2})]^{-1}\asymp n^{-1}m^{2\epsilon a+1}\asymp
n^{-2r/(2r+2a+\epsilon^{-1})}.
\]
(2) (severely ill-posed case): When $\phi(t)=\exp(-t^{-a/2})$,
$\nu_k\asymp k^\epsilon$ for some $a, \epsilon>0$, we have
\[
n^{-1}\sum_{k=1}^m\nu_k^{2r}[\phi(\nu_k^{-2})]^{-1}=n^{-1}\sum_{k=1}^m
k^{2\epsilon r}\exp(k^{a\epsilon})\asymp n^{-1}m^{2\epsilon
r}\exp(m^{a\epsilon})
\]
means that we have to choose $m=c\log(n)^{1/a\epsilon}$ with a
sufficiently small $c>0$. The resulting lower bound is
\[ \delta_n \asymp n^{-1}\sum_{k=1}^m [\phi(\nu_k^{-2})]^{-1}\asymp
n^{-1}\exp(m^{a\epsilon})\asymp m^{-2\epsilon r}\asymp (\log
n)^{-2r/a}.
\]
\end{proof}

\begin{proof}[Proof of Proposition \ref{prop-npir}]
We have
\[\E[\hat\eta_k]=\E[(Kh)(W_i)((K^\ast)^{-1}u_k)(W_i)]=\scapro{Kh}{(K^\ast)^{-1}u_k}_W
=\scapro{h}{u_k}_X
\]
and
\begin{align*}
\Var(\hat\eta_k)&=\frac1n\Var\Big((Kh)(W_i)((K^\ast)^{-1}u_k)(W_i)+V_i((K^\ast)^{-1}u_k)(W_i)\Big)\\
&\le
2n^{-1}\Big(\norm{Kh}_\infty^2\E[((K^\ast)^{-1}u_k)^2(W_i)]+\E[V_i^2]\E[((K^\ast)^{-1}u_k)^2(W_i)]\Big)\\
&\le 2n^{-1}(S^2+\sigma_1^2)\norm{(K^\ast)^{-1}u_k}_W^2.
\end{align*}
From $\norm{K g}_W\ge c\norm{[\phi(B^{-2})]^{1/2}g}_X$ for all $g\in
L_X^2$ we infer by duality $\norm{(K^\ast)^{-1} g}_W\le
c^{-1}\norm{[\phi(B^{-2})]^{-1/2}g}_X$ for all $g\in\ran(K^\ast)$.
Hence,
\[ \E_h[\norm{\hat h_n-h}_X^2]\le 2n^{-1}(S^2+\sigma_1^2)c^{-2}\sum_{k=1}^m
[\phi(\nu_k^{-2})]^{-1}+\sum_{k=m+1}^\infty \scapro{h}{u_k}_X^2.
\]
From $h\in{\cal H}(r,R)$ we have the bias estimate
\[\sum_{k=m+1}^\infty \scapro{h}{u_k}_X^2 \le \nu_{m+1}^{-2r}R^2.
\]

When choosing $m$ as for the lower bound, then the variance term
matches the lower bound in order and the estimator $\hat h_n$
attains the minimax-rate provided the bias term is not of larger
order. This is equivalent to requiring for some uniform constant
$c>0$ that $\nu_{m+1}^{2r}n^{-1}\sum_{k=1}^m
[\phi(\nu_k^{-2})]^{-1}\ge c$, which in turn follows from
$\nu_{m+1}\ge \nu_k$ for $k\le m$ and
$n^{-1}\sum_{k=1}^m\nu_k^{2r}[\phi(\nu_k^{-2})]^{-1}\asymp 1$.

(1) For mildly ill-posed case with $\phi(t)=t^a$, $\nu_k\asymp
k^\epsilon$, we have
\[
n^{-1}\sum_{k=1}^m[\phi(\nu_k^{-2})]^{-1}=n^{-1}\sum_{k=1}^m
k^{2a\epsilon}\asymp n^{-1}m^{2a\epsilon +1} \asymp m^{-2r\epsilon}
\]
by setting $m\asymp n^{1/(2\epsilon r+2\epsilon a+1)}$. Thus we
obtain the upper bound: $\delta_n \asymp m^{-2r\epsilon}\asymp
n^{-2r/(2r+2a+\epsilon^{-1})}$.

(2) For severely ill-posed case with $\phi(t)=\exp(-t^{-a/2})$,
$\nu_k\asymp k^\epsilon$, we have
\[
n^{-1}\sum_{k=1}^m [\phi(\nu_k^{-2})]^{-1}=n^{-1}\sum_{k=1}^m
\exp(k^{a\epsilon})\asymp n^{-1}\exp(m^{a\epsilon})\asymp
m^{-2r\epsilon}
\]
by setting $m=c\log(n)^{1/a\epsilon}$ with a sufficiently small
$c>0$. Thus we obtain the upper bound: $\delta_n \asymp
m^{-2r\epsilon}\asymp (\log n)^{-2r/a}$.
\end{proof}

\begin{proof}[Proof of Theorem \ref{thm-npiv}]
Given Assumption \ref{npiv-a1} ($\{\psi_k\}$ is a Riesz basis
associated with the operator $B$), there is a bounded invertible
operator $\bar{B}$ on $L^2_X$ such that $\bar{B}\psi_k = u_k$ for
all $k$. This implies that $\mathcal{H}_{m(n)} =
\spann\{u_1,...,u_{m(n)}\}$. Denote $\Pi_{m(n)} (h)$ as the
projection of $h\in \cal H (r,R)$ onto $\mathcal{H}_{m(n)}$. Then
\[
\norm{\hat h_n-h}_X^2 \le 2\{\norm{\Pi_{m(n)}(h)-h}_X^2 + \norm{\hat
h_n-\Pi_{m(n)}(h)}_X^2\}.
\]
As in Blundell, Chen and Kristensen (2007), we define $\tau _{n}$ as
a \textit{sieve measure of ill-posedness}:
\[
\tau _{n} := \sup_{h\in \mathcal{H}_{m(n)}:h\neq
0}\frac{\norm{h}_X}{\norm{Kh}_W} = \sup_{h\in
\spann\{u_1,...,u_{m(n)}\}:h\neq 0}\frac{\norm{h}_X}{\norm{Kh}_W},
\]
which is well defined under the conditions for identification. Then
\[
\norm{\hat h_n-\Pi_{m(n)}(h)}_X \le \tau_n \times \norm{K[\hat
h_n-\Pi_{m(n)}(h)]}_W.
\]
Under Assumption \ref{npiv-a2}, by the definition of $\hat{h}_n$ and
applying  Claims 2 and 3 in Blundell, Chen and Kristensen (2007), we
have:
\[
\norm{\hat h_n-\Pi_{m(n)}(h)}_X \le \tau_n \times \{
O_{p}(J_n^{-r_K} +\sqrt{(J/n)}+\norm{K[h-\Pi_{m(n)}(h)]}_W)\},
\]
where the $O_p ()$ holds uniformly over $h\in \cal H (r,R)$.

By definition of $\tau_n$ we have:
\[
{\tau}^2_{n} \le \sup_{h\in \spann\{u_1,...,u_{m(n)}\}:h\neq
0}\frac{\norm{h}^2_X}{\norm{[\phi(B^{-2})]^{1/2}h}^2_X} \le
[\phi(\nu_{m(n)}^{-2})]^{-1},
\]
where the first inequality is due to Assumption \ref{npir-a2} (the
reverse link condition), and the second inequality holds because
$\nu_k$ is increasing in $k$ and $\phi(t)$ is non-decreasing
function in $t\ge 0$.

By definition of $\tau_n$ we have under Assumptions \ref{h-smooth},
\ref{K-smooth}, \ref{npir-a2} and $\lim_{n}\frac{J_n}{m(n)}=c \in
(1, \infty)$ and $J_n > m(n)$, we obtain:
\[
{\tau}^2_{n} \norm{K[h-\Pi_{m(n)}(h)]}^2_W \le
\norm{h-\Pi_{m(n)}(h)}^2_X \le R^2 \nu_{m(n)+1}^{-2r},
\]
thus
\[
\norm{\hat h_n-h}_X^2 \le C' \max\Big\{\nu_{m(n)+1}^{-2r},~
\frac{J_n}{n} {\tau}^2_n \Big\}\le C \max\Big\{\nu_{m(n)+1}^{-2r},~
\frac{m(n)}{n} [\phi(\nu_{m(n)}^{-2})]^{-1} \Big\}
\]
uniformly over $h\in {\cal H} (r,R)$ except on an event whose
probability tends to zero as $n\uparrow \infty$.
\end{proof}
\end{document}